\documentclass[a4paper,10pt]{article}

\usepackage{amssymb}
\usepackage{graphicx}
\usepackage{graphics}

\def\yen{\hbox{iftdir\yoko\fi
\setbox0=\hbox{Y}Y\kern-.97\wd0\vbox{\hrule height.lex width.98\wd0
\kern.33ex\hrule height.lex width.98\wd0\kern.45ex}}}

\def\yen{{\setbox0=\hbox{Y}Y\kern-.97\wd0\vbox{hrule height.lex width.98%
\wd0\kern.33ex\hrule height.lex width.98\wd0\kern.45ex}}}

\def\np{\newpage}    
\begin{document}

\newtheorem{thm}{Theorem}[section]
\newtheorem{prop}{Lemma}

\def\g{\gamma}    
\def\ve{\varepsilon}   
\def\vp{\varphi}
\def\si{\mathrm{sin}} 
\def\co{\mathrm{cos}}
\def\i{\hookrightarrow}
\def\e{\hookrightarrow}
\def\l{\longrightarrow} 
\def\ttt{\longmapsto}
\def\f{\flushpar }
\def\nl{\newline }
\def\np{\newpage }
\def\x{\times } 
\def\te{^t \hskip-1mm }

\title{
Ribbon-moves of 2-knots: the Farber-Levine 
pairing and the Atiyah-Patodi-Singer- 
Casson-Gordon-Ruberman 
$\widetilde\eta$-invariants of 2-knots
}
\author{
 Eiji Ogasa
\footnote{
This research was partially supported by Research Fellowships
of the Promotion of Science for Young Scientists.
\nl{\bf Keyword:}
2-knots, ribbon-moves of 2-knots, the $\widetilde\eta$-invariants of 2-knots, 
the Farber-Levine pairing of 2-knots, the Alexander module 
}\\
 ogasa@ms.u-tokyo.ac.jp\\
Department of Physics, 
University of Tokyo\\ 
Hongo, Tokyo 113, JAPAN\\
}
\date{}
\maketitle

\baselineskip11pt
{\bf Abstract. }  
Let $K$ and $K'$ be 2-knots. 
Suppose that $K$ and $K'$ are ribbon-move equivalent. 
Then there is an isomorphism 
${\mathrm {Tor}} H_1(\widetilde X_K; {\bf Z}) \cong 
{\mathrm {Tor}} H_1(\widetilde X_{K'}; {\bf Z})$
as $ {\bf Z}[t,t^{-1}]$-modules. 
Furthermore the Farber-Levine pairing for $K$ is equivalent to that of $K'$. 

Let $K$ be a 2-knot which is ribbon-move equivalent to the trivial knot. 
Then 
the Atiyah-Patodi-Singer-Casson-Gordon-Ruberman 
$ {\bf Q}/{\bf Z}$-valued $\widetilde\eta$-invariants 
 $\tilde\eta(K, \quad)$ for ${\bf Z_d}$ is zero.
($d\in {\bf N}$ and $d\geq2$).

\section{ Ribbon-moves of 2-knots }

In this paper we discuss ribbon-moves of 2-knots.  
In this section we review the definition of ribbon-moves.

An  {\it (oriented) 2-(dimensional) knot}
 is a smooth, oriented submanifold $K$ of $S^4$ 
which is diffeomorphic to the $2$-sphere. 
We say that 2-knots $K_1$ and $K_2$ are {\it equivalent} 
if there exists an orientation preserving diffeomorphism  
$f:$ $S^4$ $\rightarrow$ $S^4$ 
such that $f(K_1)$=$K_2$  and 
that $f | _{K_1}:$ $K_1$ $\rightarrow$ $K_2$ is 
an orientation preserving diffeomorphism.    
Let $id:S^4$ $\rightarrow$ $S^4$ be the identity. 
We say that 2-knots $K_1$ and $K_2$ are {\it identical}  
if  $id(K_1)$=$K_2$ and 
that $id | _{K_1}:K_1\rightarrow K_2$ is 
an orientation preserving diffeomorphism.

\noindent{\bf Definition 1.1.}
Let $K_1$ and $K_2$ be 2-knots in $S^4$. 
We say that $K_2$ is obtained from $K_1$ by one {\it ribbon-move } 
if there is a 4-ball $B$ of $S^4$ with the following properties.  

(1) 
$K_1-(B\cap K_1)=K_2-(B\cap K_2)$. 

This diffeomorphism map is orientation preserving. 

(2) 
$B\cap K_1$ is drawn as in Figure 1.1.
$B\cap K_2$ is drawn as in Figure 1.2.

We regard $B$ as 
(a close 2-disc)$\times[0,1]\times\{t| -1\leq t\leq1\}$.
We put $B_t=$(a close 2-disc)$\times[0,1]\times\{t \}$.  
Then $B=\cup B_t$. 
In Figure 1.1 and 1.2, we draw $B_{-0.5}, B_{0}, B_{0.5}$ $\subset B$. 

We draw $K_1$ and $K_2$ by the bold line. 
The fine line denotes $\partial (B_t)$. 
  
$B\cap K_1$ (resp. $B\cap K_2$) is diffeomorphic to 
$D^2\amalg (S^1\times [0,1])$.

$B\cap K_1$ has the following properties:  
$B_t\cap K_1$ is empty for $-1\leq t<0$ and $0.5<t\leq1$.
$B_0\cap K_1$ is diffeomorphic to 
$D^2\amalg(S^1\times [0,0.3])\amalg(S^1\times [0.7,1])$. 
$B_{0.5}\cap K_1$ is diffeomorphic to $(S^1\times [0.3,0.7])$. 
$B_t\cap K_1$ is diffeomorphic to $S^1\amalg S^1$ for $0<t<0.5$.

$B\cap K_2$ has the following properties:  
$B_t\cap  K_2$ is empty for $-1\leq t<-0.5$ and $0<t\leq1$.
$B_0\cap K_2$ is diffeomorphic to 
$D^2\amalg(S^1\times [0, 0.3])\amalg(S^1\times [0.7, 1])$. 
$B_{-0.5}\cap  K_2$ is diffeomorphic to $(S^1\times [0.3, 0.7])$. 
$B_t\cap  K_2$ is diffeomorphic to $S^1\amalg S^1$ for $-0.5<t<0$. 

We do not assume which the orientation of $B\cap K_1$ (resp. $B\cap K_2$ ) is. 

\hskip3cm Figure 1.1.

\hskip3cm Figure 1.2.

Suppose that $K_2$ is obtained from $K_1$ by one ribbon-move 
and that $K'_2$ is equivalent to $K_2$.   
Then we also say that $K'_2$ is obtained from $K_1$ 
by one {\it ribbon-move}.   

If $K_1$ is obtained from $K_2$ by one ribbon-move,  
then we also say that $K_2$ is obtained from $K_1$ by one {\it ribbon-move}.

\noindent
{\bf  Definition 1.2.}
2-knots $K_1$ and $K_2$ are said to be {\it ribbon-move equivalent} 
if there are 2-knots 
$K_1=\bar{K}_1, \bar{K}_2,...,\bar{K}_{p-1},\bar{K}_p=K_2$  
 ($p\in{\bf N}, p\geq2$) such that 
$\bar{K}_i$ is obtained from $\bar{K}_{i-1}$ $(1< i\leq p)$ by one ribbon-move.

In this paper we discuss  the following problems.

\noindent
{\bf Problem 1.3. } 
Let $K_1$ and $K_2$ be 2-knots.
Consider a necessary (resp. sufficient, necessary and sufficient )
condition that $K_1$ and $K_2$ are ribbon-move equivalent.

In 
\cite{O}
 the author proved:  

\vskip2mm
\noindent
{\bf Theorem 1.4. }
(\cite{O})

{\it 
(1)If 2-knots $K$ and $K'$ are ribbon-move equivalent, then 
\vskip2mm
\hskip4cm$\mu(K)=\mu(K')$. 
\vskip2mm

(2)Let $K_1$ and $K_2$ be 2-knots in $S^4$. 
Suppose that $K_1$ are ribbon-move equivalent to $K_2$.   
Let $W_i$ be arbitrary Seifert hypersurfaces for $K_i$. 
Then the torsion part of 
$\{H_1(W_1)\oplus H_1(W_2)\}$ is congruent to 
$G\oplus G$ for a finite abelian group $G$. 

(3)Not all 2-knots are ribbon-move equivalent to the trivial 2-knot.   

(4)The inverse of (1) is not true.  The inverse of (2) is not true. 
}





 


\section{ Main results }

\noindent
{\bf Theorem 2.1.}
{\it 
Let $K$ and $K'$ be 2-knots. 
Suppose that $K$ and $K'$ are ribbon-move equivalent. 
Then we have: 

(1)There is an isomorphism 
${\mathrm {Tor}} H_1(\widetilde X_K; {\bf Z}) \cong 
{\mathrm {Tor}} H_1(\widetilde X_{K'}; {\bf  Z})$
as $  {\bf Z}[t,t^{-1}]$-modules. 

(2)The Farber-Levine pairing on 
${\mathrm{Tor}} H_1(\widetilde X_K;  {\bf Z})$ is equivalent to 
that on 
${\mathrm{Tor}} H_1(\widetilde X_{K'};  {\bf Z})$. 
}

\noindent
{\bf Theorem 2.2. }
{\it  Let $K$ be a 2-knot. 
Suppose that $K$ is ribbon-move equivalent to the trivial knot. 
Then $\tilde\eta(K, \quad)$ for ${\bf Z_d}$ is zero.
($d\in {\bf N}$ and $d\geq2$). 
}


{\bf Note.}
We review the Alexander module in \S3.
We review the Farber-Levine pairing in \S4. 
We review the Atiyah-Patodi-Singer-Casson-Gordon-Ruberman 
$  {\bf Q}/  {\bf Z}$-valued 
$\widetilde\eta$-invariants of 2-knots in \S5.

\section{  The Alexander module}  \label{the Alexander module}

See 
\cite{L66}
 \cite{R}
 for detail. 

Let $K$ be a 2-knot $\subset S^4$. 
Let $N(K)$ be the tubular neighborhood of $K$ in $S^4$. 
Let   
$\alpha:\pi_1(\overline{S^4-N(K)})\rightarrow 
H_1(\overline{S^4-N(K)};  {\bf Z})$ 
be the abelianization. 
Note that any nonzero cycle $x\in H_1(S^4-N(K);  {\bf Z})$ 
is oriented naturally 
by using the orientation of $K$ and that of $S^4$.  
We define the canonical isomorphism  
$\beta:H_1(\overline{S^4-N(K)};  {\bf Z})\rightarrow  {\bf Z}$ 
by using these orientations of $x$. 
Let $\widetilde{X}^{\infty}_K$ be 
the covering space associated with 
$\beta\circ\alpha:\pi_1(\overline{S^4-N(K)})\rightarrow  {\bf Z}$. 
We call $\widetilde{X}^{\infty}_K$ 
the {\it canonical infinite cyclic covering space}
of the complement $\overline{S^4-N(K)}$ of $K$.  
Then $H_i(\widetilde{X}^{\infty}_K;{\bf Z})$ is regarded 
as a ${\bf Z}[t,t^{-1}]$-module 
by using the covering translations
$\widetilde{X}^{\infty}_K\rightarrow\widetilde{X}^{\infty}_K$. 
This ${\bf Z}[t,t^{-1}]$-module $H_i(\widetilde{X}^{\infty}_K;{\bf Z})$ 
is called the {\it Alexander module}.

\section{ lk(\quad,\quad) for 2-knots }

See 
\cite{F}
 \cite{L77}
 for detail. 

Firstly we review of lk(\quad,\quad) for closed oriented 3-manifolds.  
Let $M$ be a closed oriented 3-manifold. 
Let $x, y\in$Tor$H_1(M; {\bf Z})$. 
Let $n$ (resp. $m$) be a natural number. 
Let $n$ (resp. $m$) be the order of $x$ (resp. $y$). 
Let $X$ (resp. $Y$) be a circle embedded in $M$ 
such that $[X]=x$ (resp. $[Y]=y$).
Let $X\cap Y=\phi$. 
Then there is an immersion map 
$f:F\longrightarrow M$ 
such that 

\hskip1cm
(1) $F$ is an oriented compact surface and $\partial F$ is one circle.  

\hskip1cm
(2) $f({\mathrm Int} F)$ is an embedding. 

\hskip1cm
(3) $f(\partial F)=X$ and deg$(f\vert_{\partial F})$=$n$. 

\hskip1cm
(4) $f({\mathrm Int} F)$ is transverse to $Y$.

Let $f(F)\cap Y=P_1\amalg...\amalg P_\alpha$. 
(Note $P_i$ is a point.) 
We give $P_i$ a signature $\varepsilon_i$ 
by using the orientation of $f(F)$, that of $Y$, and that of $M$.   

Define lk$(x,y)=\frac{1}{n}\Sigma_{i=1}^{\alpha}\varepsilon_i$ 
$\in {\bf Q}/{\bf Z}$. 

\vskip3mm
\noindent
{\bf Proposition.} {\it lk(y,x)=lk(x,y). } 
\vskip3mm

There is an immersion map $g:G\longrightarrow M$ such that 

\hskip1cm(1) $G$ is an oriented compact surface and 
$\partial G$ is one circle.

\hskip1cm(2) $g({\mathrm Int} G)$ is an embedding. 

\hskip1cm(3) $f(\partial G)=Y$ and deg$(g\vert_{\partial G})$=$m$. 

\hskip1cm(4) $f({\mathrm Int} G)$ is transverse to $X$.  

\hskip1cm(5) Int$F$ is transverse to Int$G$.

Let $g(G)\cap X=Q_1\amalg...\amalg Q_\beta$.  
(Note $Q_j$ is a point. )
We give $Q_j$ a signature $\sigma_j$ 
by using the orientation of $g(G)$, that of $X$, and that of $M$.   
Let 
$f({\mathrm Int} F)\cap g({\mathrm Int}G)=R_1\amalg...\amalg R_\gamma$.   
(Note $R_k$ is a compact open 1-manifold. )
We give $R_k$ a signature $\tau_k$ 
by using the orientation of $g(G)$, that of $f(F)$, and that of $M$.   
Then we have: 
 
\hskip1.5cmlk$(y,x)$
$=\frac{1}{m}\Sigma_{j=1}^{\beta}\sigma_j$  
$=\frac{1}{m\cdot n}\Sigma_{k=1}^{\gamma}\tau_k$ 
$=\frac{1}{n}\Sigma_{i=1}^{\alpha}\varepsilon_i$ 
$=$lk$(x,y)$.

Secondly we review lk(\quad,\quad) for 2-knots. 
Let $K$ be a 2-knot. 
Let $\widetilde{X}^{\infty}_K$ be 
the canonical infinite cyclic covering space 
of the complement $\overline{S^4-N(K)}$ of $K$. 
Let $x, y\in$ Tor$H_1(\widetilde{X}^{\infty}_K; {\bf Z})$.  
We define lk$(x, y)$.

Let $p$ be the natural projection map $\widetilde{X}^{\infty}_K\rightarrow X$. 
Let $V$ be a Seifert hypersurface for $K$. 
Let $V_\xi$ be one connected component of 
$p^{-1}(V)=\amalg^{\infty}_{-\infty}V_i$. 
The natural inclusion map 
$V_\xi\rightarrow \widetilde{X}^{\infty}_K$ induces 
the homomorphism 
$\iota: H_1(V_\xi)\rightarrow H_1(\widetilde{X}^{\infty}_K)$.  
Theorem 7.3 of 
\cite{F}
  and its proof essentially say that 
$\iota: {\mathrm {Tor}}H_1(V_\xi)\rightarrow 
{\mathrm {Tor}}H_1(\widetilde{X}^{\infty}_K)$ is onto.

Let $\hat{V_\xi}$ be the closed oriented 3-manifold 
which is obtained from $V_\xi$ 
by attaching a 3-dimensional 3-handle along $\partial V_\xi$.
The natural inclusion map $V_\xi \rightarrow \hat{V_\xi}$ induces 
 $\gamma:H_1(V_\xi; {\bf Z})\stackrel{\cong}\to 
 H_1(\hat{V_\xi};{\bf Z})$. 
Then the map 
$(\iota\circ\gamma^{-1}):  
{\mathrm{Tor}}H_1(\hat{V_\xi}; {\bf Z})\rightarrow 
{\mathrm{Tor}}H_1(\widetilde{X}^{\infty}_K; {\bf Z})$ 
is onto.  
Let $(\iota\circ\gamma^{-1})(x')=x$ and $(\iota\circ\gamma^{-1})(y')=y$.  
Define lk($x, y$) for the 2-knot $K$ to be lk$(x', y' )$ 
for the 3-manifold $\hat{V_\xi}$. 
Theorem 7.3 of 
\cite{F} 
 and its proof essentially say that  
lk($x, y$) for the 2-knot $K$ is independent of the choice of $V$.

\section{ $\widetilde\eta(\quad,\quad)$ of 2-knots }

See  $\widetilde\eta(\quad,\quad)$ for 
\cite{APS}  
 \cite{CG1} \cite{CG2}   
 \cite{Ru83}
 for detail.

We firstly define $\widetilde\eta(\quad ,\quad )$ 
of closed oriented 3-manifolds. 
(See P.571 of 
\cite{Ru83}.)      

Let $A$ be an oriented compact 4-manifold. 
For $d$ an integer, set $\omega=e^{\frac{2\pi i}{d}}$, 
and suppose 
$\phi:A\rightarrow K({\bf Z_d}, 1)$ is a continuous map. 
Then $\phi$ corresponds to a class $\phi\in H^1(A; {\bf Z_d})$=
Hom$(H_1(A; {\bf Z}); {\bf Z_d})$ and induces a cyclic cover 
$\widetilde A\rightarrow A$ with a specific choice 
$T:\widetilde A\rightarrow \widetilde A$ of a generator 
of the covering translations. 
Let $\overline{H}_2(A,\phi)$=
$\omega$-eigenspace of $T_*$ acting on $H_2(\widetilde A; {\bf C})$. 
The intersection form on $A$ induces a Hermitian pairing on 
$H_2(A; {\bf C})$; 
$<x\otimes\alpha, y\otimes\beta>$=$(x\cdot y)\otimes\alpha\overline\beta$.
Define $\overline{\sigma}( A, \phi )$=the signature of $< , >$ 
restricted to $\overline{H}_2(A,\phi)$.

Let $M$ be an oriented closed 3-manifold. 
Let $\psi\in H^1(M; {\bf Z_d})$ be a homomorphism. 
By bordism theory, 
 $n\cdot(M, \psi)=\partial(W, \phi)$ 
for a compact oriented 4-manifold $W$.  
(See e.g. 
\cite{CF} 
)
Define 
$\widetilde\eta(M,\psi)=
\frac{1}{n}(\overline\sigma(W, \phi)-\sigma(W))$.

Secondly we define $\widetilde\eta(\quad,\quad)$ of 2-knots. 
(See 
\cite{Ru83}
\cite{Ru88}
)
Let $K$ be a 2-knot.  
%
%
%
%
Take 
$\widetilde{X}^{\infty}_K$,  
$V$, $V_\xi$,  
$\iota: H_1(V_\xi;{\bf Z})\rightarrow H_1(\widetilde{X}^{\infty}_K;{\bf Z})$, 
$\hat{V_\xi}$,  
$\gamma:H_1(V_\xi; {\bf Z})\cong H_1(\hat{V_\xi}; {\bf Z})$ 
as in 
\S4. %
Let $\nu$ be a homomorphism 
$H_1(\widetilde{X}^{\infty}_K; {\bf Z})\rightarrow  {\bf Z_k}$.   

Then we have 

$$
H_1(\hat{V_\xi}; {\bf Z})\stackrel{\gamma,\cong}\longleftarrow
H_1(V_\xi; {\bf Z})\stackrel{\iota}\to
H_1(\widetilde{X}^{\infty}_K; {\bf Z})\stackrel{\nu}\rightarrow
{\bf Z_k}
$$

Then we have 
$\nu\circ\iota\circ(\gamma^{-1}):H_1(\hat{V_\xi}; {\bf Z})\to{\bf Z_k}$. 
Take 
$\widetilde\eta(\hat{V_\xi}, \nu\circ\iota\circ(\gamma^{-1}))
\in {\bf Q}$.  
Let $\pi$ be the natural projection $ {\bf Q}\to{\bf Q}/{\bf Z}$. 
\cite{Ru83}  
\cite{Ru88}  
says that 
\newline
$\pi(\widetilde\eta(\hat{V_\xi}, \nu\circ\iota\circ(\gamma^{-1})))$
is independent of the choice of $V$. 
We define 
$\widetilde\eta(K, \nu)\in{\bf Q}/{\bf Z}$ to be 
$\pi(\widetilde\eta(\hat{V_\xi}, \nu\circ\iota\circ(\gamma^{-1}))).$

\noindent
{\bf Note. } 
\cite{CO}  3
\cite{L94} 
\cite{GL}  
etc. also 
apply the G-signature theorem and the $\eta$-invariants 
to $n$-knots ($n=1,2,3,...$).

\section{ Proof of main results}


Let $X$ be $\overline{S^4-N(K)}$. 
Let $V$ be a Seifert hypersurface for $K$. 
Suppose $V$ also denote $\overline{V-N(K)}$. 
Let $Y=X-V$.
Let $V\x[-1,1]$ be the tubular neighborhood of $V$ in $X$. 
Suppose $Y$ also denote $\overline{X-(V\x[-1,1])}$. 

Let $\widetilde{X}^{\infty}_K$ be 
the canonical infinite cyclic covering space.
There is the natural projection map $p:\widetilde{X}^{\infty}_K\rightarrow X$. 
Let $p^{-1}(Y)=\amalg^{\infty}_{-\infty}Y_i$. 
Let $p^{-1}(V)=\amalg^{\infty}_{-\infty}V_i$. 
Let $p^{-1}(V\x[-1,1])=\amalg^{\infty}_{-\infty}(V_i\x[-1,1])$. 
Let $\partial Y_i=(V_{i-1}\x\{1\})\amalg (V_i\x\{-1\})$.

\np
There is the Meyer-Vietoris exact sequences
of {\bf Z}-homology groups:


$$H_i(V\x\{-1,1\};  {\bf Z}) $$
$$\downarrow^{f_i}$$
$$H_i(V\x[-1,1];  {\bf Z})\oplus H_i(Y;  {\bf Z})  $$
$$\downarrow^{g_i}$$
$$H_i(X;  {\bf Z})$$

and 


$$\oplus^{\infty}_{j=-\infty} H_i(V_j\x\{-1,1\};  {\bf Z})$$
$$\downarrow^{\widetilde {f_i}} $$
$$ (\oplus^{\infty}_{j=-\infty} H_i(V_j\x[-1,1] ;  {\bf Z}))
\oplus (\oplus^{\infty}_{j=-\infty}H_i(Y_j  ;  {\bf Z}))$$
$$\downarrow^{\widetilde {g_i}} $$
$$H_i(\widetilde{X}^{\infty}_K ;  {\bf Z})$$

Furthermore the second one is regarded as an exact sequence 
of ${\bf Z}[t,t^{-1}]$-modules by using the covering tranlations.

\noindent
{\bf Claim 6.1. } 
{\it 
There is an exact sequence of 
${\bf Z}[t,t^{-1}]$-modules:}

$$\oplus^{\infty}_{j=-\infty}\mathrm{Tor}H_1(V_j\x\{-1,1\}; {\bf  Z})$$
$$\downarrow^{\widetilde {f^{\mathrm{tor}}}}$$
$$ (\oplus^{\infty}_{j=-\infty}\mathrm{Tor}H_1(V_j\x[-1,1]; {\bf Z})) 
\oplus (\oplus^{\infty}_{j=-\infty}\mathrm{Tor}H_1(Y_j; {\bf Z}))$$
$$\downarrow^{\widetilde {g^{\mathrm{tor}}}} $$
$$\mathrm{Tor}H_1(\widetilde{X}^{\infty}_K;{\bf Z})$$
$$\downarrow $$
$$0$$

{\it 
where $\widetilde {f^{\mathrm{tor}}}=$
$\widetilde {f_1}\vert_
{\oplus^{\infty}_{j=-\infty}\mathrm{Tor}H_1(V_j\x\{-1,1\}; {\bf  Z})}$ 

and  $\widetilde {g^{\mathrm{tor}}}=$
$\widetilde {g_1}\vert_
{(\oplus^{\infty}_{j=-\infty}\mathrm{Tor}H_1(V_j\x[-1,1]; {\bf Z})) 
\oplus (\oplus^{\infty}_{j=-\infty}\mathrm{Tor}H_1(Y_j; {\bf Z}))}$.  
}

\noindent{\bf Note. }
P.765 of 
\cite{F} 
says that 
\cite{K}   
proved 
$\vert{\mathrm {Tor}}H_1(\widetilde{X}^{\infty}_K; {\bf Z})\vert<\infty$.

\noindent{\bf Proof.} 
By using the covering translations, we regard 
$\oplus^{\infty}_{j=-\infty}\mathrm{Tor}H_1(V_j\x\{-1,1\}; {\bf  Z})$, 
$ (\oplus^{\infty}_{j=-\infty}\mathrm{Tor}H_1(V_j\x[-1,1]; {\bf Z})) 
\oplus (\oplus^{\infty}_{j=-\infty}\mathrm{Tor}H_1(Y_j; {\bf Z}))$, and 
$\mathrm{Tor}H_1(\widetilde{X}^{\infty}_K;{\bf Z})$ 
as ${\bf Z}[t,t^{-1}]$-modules. 
 Furthermore we regard 
$\widetilde {f^{\mathrm{tor}}}$ and 
$\widetilde {g^{\mathrm{tor}}}$
as homomorphisms of ${\bf Z}[t,t^{-1}]$-modules.

Take $V_\xi$, $\iota$ as in \S4. 
Theorem 7.3 of 
\cite{F}  
and its proof essentially say that 
$\iota: {\mathrm {Tor}} H_1(V_\xi)\rightarrow 
{\mathrm {Tor}} H_1(\widetilde{X}^{\infty}_K)$ 
 is onto.  
Note 
$\widetilde {g^{\mathrm {tor}}}\vert_
{{\mathrm {Tor}}H_1(V_\xi\times[-1,1])}=\iota$. 
Hence $\widetilde {g^{\mathrm {tor}}}$ is onto. 
Therefore


$$(\oplus^{\infty}_{j=-\infty}\mathrm{Tor}H_1(V_j\x[-1,1]; {\bf Z})) 
\oplus (\oplus^{\infty}_{j=-\infty}\mathrm{Tor}H_1(Y_j; {\bf Z}))$$
$$\downarrow^{\widetilde {g^{\mathrm{tor}}}} $$
$$\mathrm{Tor}H_1(\widetilde{X}^{\infty}_K; {\bf Z})$$
$$\downarrow $$
$$0$$

is exact. 

Let 
$x\in(\oplus^{\infty}_{j=-\infty}\mathrm{Tor}H_1(V_j\x[-1,1]; {\bf Z})) 
\oplus (\oplus^{\infty}_{j=-\infty}\mathrm{Tor}H_1(Y_j; {\bf Z}))$
such that $\widetilde {g^{\mathrm{tor}}}(x)=0$. 
Then there is 
$y\in\oplus^{\infty}_{j=-\infty}\mathrm{Tor}H_1(V_j\x\{-1,1\};{\bf Z})$
such that ${\widetilde {f_1}(y)}=x$. 
Let $n$ be the order of $x$. 
Then ${\widetilde {f_1}}(n\cdot y)=n\cdot x=0$. 
We prove that 
$y\in\oplus^{\infty}_{j=-\infty}\mathrm{Tor}H_1(V_j\x\{-1,1\}; {\bf  Z})$ 
 in the following paragraphs.

There is the Meyer-Vietoris exact sequences of 
{\bf Q}-homology groups: 


(1)

$$H_i(V\x\{-1,1\};  {\bf Q}) $$
$$\downarrow^{f^Q_i}$$
$$H_i(V\x[-1,1];  {\bf Q})\oplus H_i(Y;  {\bf Q})  $$
$$\downarrow^{g^Q_i}$$
$$H_i(X;  {\bf Q})$$

and 


(2)

$$\oplus^{\infty}_{j=-\infty} H_i(V_j\x\{-1,1\};  {\bf Q})$$
$$\downarrow^{\widetilde {f^Q_i}} $$
$$ (\oplus^{\infty}_{j=-\infty} H_i(V_j\x[-1,1] ;  {\bf Q}))
\oplus (\oplus^{\infty}_{j=-\infty}H_i(Y_j  ;  {\bf Q}))$$
$$\downarrow^{\widetilde {g^Q_i}} $$
$$H_i(\widetilde{X}^{\infty}_K ;  {\bf Q})$$

By using 
$H_i(X;{\bf Q})\cong H_i(S^1;{\bf Q})$ and the sequence (1), 

$f^Q_1: H_1(V\x\{-1,1\};  {\bf Q}) 
\rightarrow 
H_1(V\x[-1,1];{\bf Q})\oplus H_1(Y;  {\bf Q})$ 
is isomorphism.

Let 
$\pi_V:H_1(V\x[-1,1];{\bf Q})\oplus H_1(Y;  {\bf Q})  
\to H_1(V\x[-1,1];{\bf Q})$ be the natural projection.  
Let 
$\pi_Y:H_1(V\x[-1,1];{\bf Q})\oplus H_1(Y;  {\bf Q})  
\to H_1(Y;  {\bf Q})$ be the natural projection.  
Suppose that 
the identity matrix $E$ represents

$\pi_V\circ \{f^Q_1\vert_{H_1(V\x\{1\};{\bf Q})}\}:$
$H_1(V\x\{1\};{\bf Q})\rightarrow H_1(V\x[-1,1];{\bf Q})$, 

the identity matrix $E$ represents 

$\pi_V\circ \{f^Q_1\vert_{H_1(V\x\{-1\};{\bf Q})}\}:$
$H_1(V\x\{-1\};{\bf Q})\rightarrow H_1(V\x[-1,1];{\bf Q})$,  

a matrix $A$ represents 

$\pi_Y\circ \{f^Q_1\vert_{H_1(V\x\{1\};{\bf Q})}\}:$
$H_1(V\x\{1\};{\bf Q})\rightarrow H_1(Y;{\bf Q})$, and 

a matrix $B$ represents 

$\pi_Y\circ \{f^Q_1\vert_{H_1(V\x\{-1\};{\bf Q})}\}:$
$H_1(V\x\{-1\};{\bf Q})\rightarrow H_1(Y;{\bf Q})$.

Then 
$f^Q_1: H_i(V\x\{-1,1\};  {\bf Q}) \rightarrow 
H_i(V\x[-1,1];{\bf Q})\oplus H_i(Y;  {\bf Q})$ is represented by

$P=\left(
\begin{array}{cc}
E&A\\
E&B\\
\end{array}
\right). $

Since $H_i(X;{\bf Q})\cong H_i(S^1;{\bf Q})$, det$P\neq0$. 

We regard 
$ H_i(V_j\x\{-1,1\};{\bf Q})$,   
$\oplus^{\infty}_{j=-\infty} H_i(V_j\x[-1,1];{\bf Q})$, and 
$\oplus^{\infty}_{j=-\infty}H_i(Y_j;{\bf Q})$ 
as ${\bf Q}[t,t^{-1}]$-modules 
by using the covering translations. 
Then 
$\widetilde {f^Q_1}$ is represented by 

$P(t)=\left(
\begin{array}{cc}
E&t\cdot A\\
E&B\\
\end{array}
\right). $

Note $P(1)=P$.
Then det$P(1)=$det$P\neq0$.
Hence det$P(t)\neq0$.
Hence 
$\widetilde {f^Q_1}$ is injective. 
Hence we have: 

if an element 
$z$ 
$\in\oplus^{\infty}_{j=-\infty}H_1(V_j\x\{-1,1\}; {\bf  Z})$  

generates 
${\bf Z}$
$\subset\oplus^{\infty}_{j=-\infty}H_1(V_j\x\{-1,1\}; {\bf  Z})$, 
$\widetilde {f_1}(z)\neq0$.  
Hence we have: 

if 
$y$ 
$\in\oplus^{\infty}_{j=-\infty}H_1(V_j\x\{-1,1\}; {\bf  Z})$ 

generates 
 ${\bf Z}$
$\subset\oplus^{\infty}_{j=-\infty}H_1(V_j\x\{-1,1\}; {\bf  Z})$, 
$\widetilde {f_1}(n\cdot y)\neq0$.  
Therefore 
\newline
$y\in\oplus^{\infty}_{j=-\infty}\mathrm{Tor}H_1(V_j\x\{-1,1\};{\bf Z})$. 
Therefore


$$\oplus^{\infty}_{j=-\infty}\mathrm{Tor}\{H_1(V_j\x\{-1,1\}; {\bf  Z})\}$$
$$\downarrow^{\widetilde {f^{\mathrm{tor}}}}$$
$$ (\oplus^{\infty}_{j=-\infty}\mathrm{Tor}\{H_1(V_j\x[-1,1]; {\bf Z})\}) 
\oplus (\oplus^{\infty}_{j=-\infty}\mathrm{Tor}\{H_1(Y_j; {\bf Z})\})$$
$$\downarrow^{\widetilde {g^{\mathrm{tor}}}} $$
$$\mathrm{Tor}\{H_1(\widetilde{X}^{\infty}_K; {\bf Z})\}$$

is exact. 
This completes the proof of Claim 6.1.

In order to prove our main theorems, 
 we use the (1,2)-pass-moves for 2-knots. 
See [13] for the (1,2)-pass-moves for 2-knots for detail. 

\vskip3mm
\noindent
{\bf Definition 6.2. } %
Let $K_1$ and $K_2$ be 2-knots in $S^4$. 
We say that $K_2$ is obtained from $K_1$ by one {\it (1,2)-pass-move } 
if there is a 4-ball $B$ $\subset S^4$ with the following properties.  
We draw $B$ as in Definition 1.1.

(1) 
$K_1-(B\cap K_1)$=$K_2-(B\cap K_2)$. 

This diffeomorphism map is orientation preserving.

(2)
$B\cap K_1$  is drawn as in 
Figure 6.1. %
$B\cap K_2$ is drawn as in 
Figure 6.2.  %

\hskip3cm Figure 6.1.   %

\hskip3cm Figure 6.2.   %

The orientation of the two discs in the 
Figure 6.1 (resp. Figure 6.2) %
is compatible with 
the orientation which is determined naturally 
by the $(x,y)$-arrows in the Figure. 
We do not assume which the orientations of the annuli in the Figures are. 

Suppose that $K_2$ is obtained from $K_1$ by one (1,2)-pass-move 
and that $K'_2$ is equivalent to $K_2$.   
Then we also say that $K'_2$ is obtained from $K_1$ by one 
{\it (1,2)-pass-move }.

If $K_1$ is obtained from $K_2$ by one (1,2)-pass-move,  then 
we also say that $K_2$ is obtained from $K_1$ by one {\it (1,2)-pass-move }.  

2-knots $K_1$ and $K_2$ are said to be {\it (1,2)-pass-move equivalent} 
if there are 2-knots 
$K_1=\bar{K}_1, \bar{K}_2,...,\bar{K}_{p-1},\bar{K}_p=K_2$  
$( p\in {\bf N}, p\geq2 )$ such that 
$\bar{K}_i$ is obtained from $\bar{K}_{i-1}$ $(1< i\leq p)$ 
by one (1,2)-pass-move. 

In 
\cite{O} 
we proved: 

\vskip3mm
\noindent
{\bf Theorem 6.3 } 
(\cite{O}) 
       {\it 
Let $K$ and $K'$ be 2-knots.  
The following conditions (1) and (2) are equivalent. 

(1)
$K$ is  (1,2)-pass-move equivalent to $K'$.

(2)
$K$ is  ribbon-move equivalent to $K'$. 
}
\vskip3mm

Let $K_<$ and $K_>$ be 2-knots. 
Suppose $K_<$ is ribbon-move equivalent to $K_>$. 
By 
Theorem 6.3, 
$K_<$ is (1,2)-pass-move equivalent to $K_>$.  

In order to prove 
Theorem 2.1 %
it suffices to prove 
Theorem 2.1 
when 
$K_<$ is obtained from $K_>$ 
by one (1,2)-pass-move in a 4-ball $B \subset S^4$.

\vskip3mm
\noindent
{\bf Claim.} {\it 
There are Seifert hypersurfaces $V_>$ for $K_>$ and $V_<$ for $K_<$ 
such that: 


(1) 
$V_>-(B\cap V_>)=V_<-(B\cap V_<)$. 

This diffeomorphism map is orientation preserving. 

(2) 
$B\cap V_>$ is drawn as in Figure 6.3.  
$B\cap V_<$ is drawn as in Figure 6.4.  
}
\vskip3mm

\noindent
{\bf Note.}  
We draw $B$ as in Definition 1.1. 
We draw $V_>$ and $V_<$ by the bold line. 
The fine line means $\partial B$. 
  
$B\cap V_>$ (resp. $B\cap V_<$) is diffeomorphic to 
$(D^2\times [2,3])\amalg (D^2\times [0,1])$. 
We can regard $(D^2\times [0,1])$ as a 3-dimensional 1-handle 
which is attached to $\partial B$. 
We can regard $(D^2\times [2,3])$ as a 3-dimensional 2-handle 
which is attached to $\partial B$.

$B\cap V_>$ has the following properties:  
$B_t\cap V_>$ is empty for $-1\leq t<0$ and $0.5<t\leq1$.
$B_0\cap V_>$ is diffeomorphic to 
$(D^2\times [2,3])\amalg(D^2\times [0,0.3])\amalg(D^2\times [0.7,1])$. 
$B_{0.5}\cap K_1$ is diffeomorphic to $(D^2\times [0.3,0.7])$. 
$B_t\cap V_>$ is diffeomorphic to $D^2\amalg D^2$ for $0<t<0.5$.

$B\cap V_<$ has the following properties: . 
$B_t\cap  V_<$ is empty for $-1\leq t<-0.5$ and $0<t\leq1$.
$B_0\cap V_<$ is diffeomorphic to 
$(D^2\times[2,3]\amalg(D^2\times [0, 0.3])\amalg(D^2\times [0.7, 1])$. 
$B_{-0.5}\cap  V_<$ is diffeomorphic to $(D^2 \times [0.3, 0.7])$. 
$B_t\cap V_<$ is diffeomorphic to $D^2\amalg D^2$ for $-0.5<t<0$. 


\hskip3cm Figure 6.3.   

\hskip3cm Figure 6.4.  

\noindent
{\bf Proof of Claim. }  
Put  $P=($the 3-manifolds in Figure 6.3$)\cap(\partial B)$. 
Note $P=($the 3-manifolds in Figure 6.4$)\cap(\partial B)$. 
Put $Q=K_>\cap(S^4-Int B^4)$. 
Note $Q=K_<\cap(S^4-Int B^4)$. 
By applying the following  Proposition to ($P\cup Q$) 
and $(S^4-Int B^4)$ the above Claim holds.

By using the obstruction theory, 
we have the following proposition. 
 ( We can prove it by applying \S III of 
\cite{S}.    
We can also prove it by generalizing Theorem 2,3 in P.49,50 of 
\cite{Ki}. )  
The author gives a proof in the Appendix. 

\vskip3mm
\noindent
{\bf Proposition.} {\it 
Let $X$ be an oriented compact $(m+2)$-dimensional manifold. 
Let $\partial X\neq\phi$. 
Let $M$ be an oriented closed $m$-dimensional manifold 
which is embedded in $X$. 
Let $M\cap\partial X\neq\phi$. 
Let $[M]=0\in H_m(X; {\bf Z})$. 
Then there is an oriented compact $(m+1)$-dimensional manifold $P$ 
such that $P$ is embedded in $X$ and that $\partial P=X$. 
}
\vskip3mm


Let $V_\asymp$ be a compact 3-manifold embedded in $S^4$ 
whose boundary is $S^2\amalg S^2$ 
with the following properties. 

(1) 
$V_\asymp-(B\cap V_\asymp)
=V_>-(B\cap V_>)=V_<-(B\cap V_<)$. 

(2) $B\cap V_\asymp\subset B_0.$    
We draw $B\cap V_\asymp$ as in Figure 6.5.

\hskip3cm Figure 6.5.

Let $N(\partial V_\asymp)$ be the tubular neighborhood of 
$\partial V_\asymp$ in $S^4$.  
Let $X=\overline{S^4-N(\partial V_\asymp)}$. 
Let $V_\asymp$ also denote $V_\asymp \cap X$. 
Let $V_\asymp\times[-1,1]$ be the tubular neighborhood of 
$V_\asymp$ in $X$. 
Let $Y_\asymp$=$\overline{X-(V_\asymp\times[-1,1])}$.



Suppose the following maps are inclusion. 
( $\varepsilon=-1,1$. ) 
The following two diagrams are commutative.

\vskip3mm
$
\begin{array}{ccc}
V_<\x\{\varepsilon\}& 
\stackrel{}\longleftarrow&  
V_\asymp\x\{\varepsilon\}\\
\downarrow ^{}& &\downarrow^{}\\     
Y_<&\stackrel{ }\longrightarrow &Y_\asymp\\
\end{array}
$

\vskip3mm\hskip28mm
$
\begin{array}{ccc}
V_\asymp\x\{\varepsilon\}& 
\stackrel{}\longrightarrow & 
V_>\x\{\varepsilon\}\\
\downarrow^{}& &\downarrow^{\alpha_>} \\     
Y_\asymp&\stackrel{}\longleftarrow & Y_>\\
\end{array}
$

\vskip3mm

By the definition of $V_\asymp$, 
$H_1(X;{\bf Z})={\bf Z}\oplus{\bf Z}$.
Let 
$f:\pi_1X\to{\bf Z}$ be an epimorphism. 
Then $f$ induces a smooth map 
$\bar f:X\to S^1.$ 
Let $q\in S^1$ be a regular value of $\bar f$. 
Suppose $\bar f^{-1}(q)=V_\asymp.$
Take the covering space $\widetilde X$ of $X$ 
associated with $\bar f$. 
Let ${\widetilde f}:\widetilde X\to X$ be the projection map. 
Put  $\amalg_{j=-\infty}^{\infty}V_{\asymp,j}
={\widetilde f}^{-1}(V_\asymp)$. 
Put  $\amalg_{j=-\infty}^{\infty}Y_{\asymp,j}
={\widetilde f}^{-1}(Y_\asymp)$. 
Let 
$\amalg_{j=-\infty}^{\infty}(V_{\asymp,j}\times[-1,1])
=\bar f^{-1}(V_\asymp\times[-1,1])$. 
Suppose 
$\partial Y_{\asymp,j}= 
(V_{\asymp,j-1}\times\{1\})\amalg (V_{\asymp,j}\times\{-1\})$.

Take $V_{>,j}$,  $V_{<,j}$, $\widetilde X^{\infty}_{K_>}$, 
$\widetilde X^{\infty}_{K_<}$ as in \S4. 
Suppose the following maps are inclusion. 
By the above commutative diagrams,  
the following diagrams are commutative.

\vskip3mm
$\begin{array}{ccc}
V_{<,j}\x\{\varepsilon\}& 
\stackrel{}\longleftarrow&  
V_{\asymp,j}\x\{\varepsilon\}\\
\downarrow ^{}& &\downarrow^{}\\     
Y_{<,j}&\stackrel{ }\longrightarrow &Y_{\asymp,j}\\
\end{array}$

\vskip3mm\hskip28mm
$\begin{array}{ccc}
V_{\asymp,j}\x\{\varepsilon\}& 
\stackrel{}\longrightarrow & 
V_{>,j}\x\{\varepsilon\}\\
\downarrow^{}& &\downarrow^{\alpha_>} \\     
Y_{\asymp,j}&\stackrel{}\longleftarrow & Y_{>,j}\\
\end{array}$

\vskip3mm

By the definition of $V_>$, $V_<$, and $V_\asymp$, it holds that 
$V_>$ (resp. $V_<$) is obtained from $V_\asymp$ by attaching one 1-handle. 
By the definition of $Y_>$, $Y_<$, and $Y_\asymp$, it holds that 
$Y_\asymp$ is obtained from $Y_>$ (resp. $Y_<$) by attaching one 3-handle. 
Hence the above commutative diagrams induce 
the following commutative diagram ($\varepsilon=-1,1$).

\[\begin{array}{ccccc}
{\mathrm {Tor}} H_1(V_{<,j}\x\{\varepsilon\} ; Z) & 
\stackrel{\beta_{<,j}, \cong}\longleftarrow&  
{\mathrm {Tor}} H_1(V_{\asymp,j}\x\{\varepsilon\} ; Z) & 
\stackrel{\beta_{>,j}, \cong} \longrightarrow & 
{\mathrm {Tor}} H_1(V_{>,j}\x\{\varepsilon\} ; Z) \\
\downarrow ^{\alpha_{<,j}}& &
\downarrow^{\alpha_{\asymp,j}}& &\downarrow^{\alpha_{>,j}} \\     
{\mathrm {Tor}} H_1(Y_{<,j} ; Z) & 
\stackrel{\gamma_{<,j}, \cong}\longrightarrow & 
{\mathrm {Tor}} H_1(Y_{\asymp,j}; Z) & 
\stackrel{\gamma_{>,j}, \cong}\longleftarrow & 
{\mathrm {Tor}} H_1(Y_{>,j}; Z) \\
\end{array}\]



\vskip3mm
The above commutative diagram induces 
the following commutative diagram of {\bf Z}-homology groups:
\vskip3mm

\noindent
{\scriptsize 
$
\begin{array}{ccccc}
\oplus_{-\infty}^\infty{\mathrm {Tor}}H_1(V_{<,j}\x\{\ -1,1\}) & 
\stackrel{\beta_{<,j}, \cong}\longleftarrow &
\oplus_{-\infty}^\infty{\mathrm {Tor}}H_1(V_{\asymp,j}\x\{\ -1,1\})   &
\stackrel{\beta_{>,j}, \cong}\longrightarrow &
\oplus_{-\infty}^\infty{\mathrm {Tor}}H_1(V_{>,j}\x\{\ -1,1\}) \\
\downarrow^{f_<} & &\downarrow^{f_\asymp}& &\downarrow^{f_>}\\
\left\{
\begin{array}{c}
\oplus_{-\infty}^\infty{\mathrm {Tor}}H_1(V_{<,j}\x[-1,1] )\\
\oplus\\
\oplus_{-\infty}^\infty{\mathrm {Tor}}H_1(Y_{<,j} ) \\
\end{array}
\right\}
& 
\stackrel{\gamma_{<,j}, \cong}\longrightarrow &
\left\{
\begin{array}{c}
\oplus_{-\infty}^\infty{\mathrm {Tor}}H_1(V_{\asymp,j}\x[-1,1] )\\
\oplus\\
\oplus_{-\infty}^\infty{\mathrm {Tor}}H_1(Y_{\asymp,j} ) \\
\end{array}
\right\}
& 
\stackrel{\gamma_{>,j}, \cong}\longleftarrow 
&
\left\{
\begin{array}{c}
\oplus_{-\infty}^\infty{\mathrm {Tor}}H_1(V_{>,j}\x[-1,1] )\\
\oplus\\
\oplus_{-\infty}^\infty{\mathrm {Tor}}H_1(Y_{>,j} ) \\
\end{array}
\right\}
\end{array}
$
}


\vskip3mm
By using the covering tlanslations, the above commutative diagram
is regarded as that of ${\bf Z}[t,t^{-1}]$-modules.  
By Claim 6.3, 
we have exact sequences of ${\bf Z}[t,t^{-1}]$-modules.  
\vskip3mm



$
\begin{array}{ccc}
\oplus_{-\infty}^\infty{\mathrm {Tor}}H_1(V_{<,j}\x\{\ -1,1\}) & 
&
\oplus_{-\infty}^\infty{\mathrm {Tor}}H_1(V_{>,j}\x\{\ -1,1\}) \\
\downarrow^{f_<} & &\downarrow^{f_>}\\
\oplus_{-\infty}^\infty
\left(
\begin{array}{c}
\mathrm {Tor}H_1(V_{<,j}\x[-1,1] )\\
\oplus\\
\mathrm {Tor}H_1(Y_{<,j} ) 
\end{array}
\right)
&
\mathrm{ and }
&
\oplus_{-\infty}^\infty
\left(
\begin{array}{c}
\mathrm {Tor}H_1(V_{>,j}\x[-1,1] )\\
\oplus\\
\mathrm{Tor}H_1(Y_{>,j} )\\
\end{array}
\right)
\\
\downarrow^{g_>}& &\downarrow^{g_>}\\
{\mathrm {Tor}}H_1(\widetilde X^{\infty}_{K_<})& 
&
{\mathrm {Tor}}H_1(\widetilde X^{\infty}_{K_>}) \\
\downarrow & &\downarrow \\
0&  & 0 \\
\end{array}
$
\vskip3mm

Therefore we have 
the following commutative diagram of ${\bf Z}[t,t^{-1}]$-modules.


\vskip3mm
$
\begin{array}{ccc}
\oplus_{-\infty}^\infty{\mathrm {Tor}}H_1(V_{<,j}\x\{\ -1,1\}) & 
\stackrel{a, \cong}\longrightarrow &
\oplus_{-\infty}^\infty{\mathrm {Tor}}H_1(V_{>,j}\x\{\ -1,1\}) \\
\downarrow^{f_<} & &\downarrow^{f_>}\\
\oplus_{-\infty}^\infty
\left(
\begin{array}{c}
\mathrm {Tor}H_1(V_{<,j}\x[-1,1] )\\
\oplus\\
\mathrm {Tor}H_1(Y_{<,j} ) 
\end{array}
\right)
&
\stackrel{b, \cong}\longrightarrow &
\oplus_{-\infty}^\infty
\left(
\begin{array}{c}
\mathrm {Tor}H_1(V_{>,j}\x[-1,1] )\\
\oplus\\
\mathrm{Tor}H_1(Y_{>,j} )\\
\end{array}
\right)
\\
\downarrow^{g_>}& &\downarrow^{g_>}\\
{\mathrm {Tor}}H_1(\widetilde X^{\infty}_{K_<})& 
\stackrel{c, \cong}\longrightarrow &
{\mathrm {Tor}}H_1(\widetilde X^{\infty}_{K_>}) \\
\downarrow & &\downarrow \\
0&  & 0 \\
\end{array}
$

\vskip3mm


In particular, 
${\mathrm {Tor}}H_1(\widetilde X^{\infty}_{K_<})
\cong 
{\mathrm {Tor}}H_1(\widetilde X^{\infty}_{K_>})$
as $ {\bf Z}[t,t^{-1}]$-module. 
This completes the proof of Theorem 2.1.(1).

Take $V_{<,\xi}$ $V_{>,\xi}$ as in \S4.  
Then there is an orientation preserving diffeomorphism 
 $h:V_{<,\xi}\to V_{>,\xi}$ with the following properties.
 
 (1) The following diagram is commutative

\[\begin{array}{ccccc}
V_{<,\xi}&&\stackrel{h}\rightarrow&&V_{>,\xi}\\
&\nwarrow &&\nearrow &\\
&&V_{\asymp,\xi}&&\\
\end{array}\] 

\hskip7mm
where $\nwarrow$ and  $\nearrow$are inclusion.

(2) $h_*=b\vert_{{\mathrm {Tor}} H_1(V_{<,\xi}; {\bf Z})}: 
{\mathrm {Tor}} H_1(V_{<,\xi}; {\bf Z}) 
\stackrel{\cong}\to  
{\mathrm {Tor}} H_1(V_{>,\xi}; {\bf Z}). $

Theorem 7.3 of 
[Farber ] 
and its proof essentially say that 

${\mathrm {Tor}} H_1(V_{<, \xi}; {\bf Z}) \stackrel{\iota_<}\to 
{\mathrm {Tor}} H_1(\widetilde{X}^{\infty}_{K_<}; {\bf Z}) $ 
and  
${\mathrm {Tor}} H_1(V_{>, \xi}; {\bf Z}) \stackrel{\iota_>}\to 
 {\mathrm {Tor}} H_1(\widetilde{X}^{\infty}_{K_>}; {\bf Z})$ 
are onto. 
Hence there are the following commutative diagram:

\[\begin{array}{ccc}
{\mathrm {Tor}} H_1(V_{<,\xi}; {\bf Z}) & 
\stackrel{h_*, \cong}\to  &
{\mathrm {Tor}} H_1(V_{>,\xi}; {\bf Z}) 
\\ 
\downarrow ^{\iota_<, \cong}&&\downarrow^{\iota_>, \cong} 
 \\ 
{\mathrm {Tor}} H_1(\widetilde{X}^{\infty}_{K_<}; {\bf Z}) & 
\stackrel{c, \cong}\to &
 {\mathrm {Tor}} H_1(\widetilde{X}^{\infty}_{K_>}; {\bf Z})
 \\ 
\end{array}\]

We can define lk(\quad,\quad) of $K_\sharp$
by using $V_{\sharp,\xi}$ ($\sharp=<,>$).

Let $!=a,b.$
Let $x_!\in{\mathrm {Tor}} H_1(\widetilde{X}^{\infty}_{K_<}; {\bf Z})$.
Then there is $x'_!\in{\mathrm {Tor}} H_1(V_{<,\xi}; {\bf Z})$ 
such that $\iota_<(x_!')=x_!$. 
Then $\iota_>\circ h(x_!')=c(x_!')$. 
Therefore 

lk$(x_a,x_b)$ for $K_<$

=lk$(x'_a,x'_b)$ for $\hat V_{<,\xi}$  
($\hat V_{<,\xi}$ is defined for $V_{<,\xi}$ as in \S4.)

=lk$(h(x'_a), h(x'_b))$ for $\hat V_{>,\xi}$
($\hat V_{>,\xi}$ is defined for $V_{>,\xi}$ as in \S4.)

=lk$(\iota_>\circ h(x'_a), \iota_>\circ h(x'_b))$ for $K_>$

=lk$(c(x'_a), c(x'_b))$ for $K_>$

Therefore the Farber-Levine pairing for $K_>$ 
coincides with that for $K_<$.  
This completes the proof of 
Theorem 2.1.(2). %


We next prove 
Theorem 2.2. %
Let $\nu$ be a homomorphism 
$H_1(\widetilde{X}^{\infty}_K; {\bf Z})\rightarrow{\bf Z_d}$. 
We consider 
$\widetilde\eta(K,\nu)$. 
Suppose $K$ is ribbon-move equivalent to the trivial 2-knot. 
Take $\iota$ and $\gamma$ as in \S5. 
Then we have the following commutative diagram.

$$
\begin{array}{ccc}
H_1(\widetilde{X}^{\infty}_K;{\bf Z})&\stackrel{\nu}\rightarrow&{\bf Z_d}\\ 
\uparrow_{\iota\circ(\gamma^{-1})} &\nearrow &\\
H_1(\hat{V_\xi};{\bf Z})&&
\end{array}
$$

By Theorem 2.1.(1), 
$H_1(\widetilde{X}^{\infty}_K;{\bf Z})\cong
{\bf Z}\oplus...\oplus {\bf Z}$.
Hence

$$
\begin{array}{ccc}
{\bf Z}\oplus...\oplus {\bf Z}&\stackrel{\nu}\rightarrow&{\bf Z_d}\\ 
\uparrow_{\iota\circ(\gamma^{-1})} &\nearrow &\\
H_1(\hat{V_\xi};{\bf Z})&&
\end{array}
$$

By an elementary discussion on homomorphisms, we have: 

$$
\begin{array}{ccc}
{\bf Z}&\stackrel{\nu'}\rightarrow&{\bf Z_d}\\ 
\uparrow_{\zeta} &\nearrow &\\
H_1(\hat{V_\xi};{\bf Z})&&
\end{array}
$$

The above homomorphism 
$H_1(\hat{V_\xi};{\bf Z})\to{\bf Z}$ is called $\zeta$.
Then \newline\noindent$\nu'\circ\zeta=$ $\nu\circ\iota\circ(\gamma^{-1})$
Then we can regard  

$\zeta\in$Hom$(H_1(\hat{V_\xi};{\bf Z}),{\bf Z})\cong 
H^1(\hat{V_\xi};{\bf Z})$

$\hskip2cm\cong  
\{\mathrm{homotopy\hskip1mm classes\hskip1mm of\hskip1mm maps}\hskip1mm  
\hat{V_\xi}\to K({\bf Z}, 1) \}$ 

$\nu'\circ\zeta\in$Hom$(H_1(\hat{V_\xi};{\bf Z}),{\bf Z_d})\cong 
H^1(\hat{V_\xi};{\bf Z_d})$

$\hskip2cm\cong
\{\mathrm{homotopy\hskip1mm classes\hskip1mm of\hskip1mm maps}
\hskip1mm    
\hat{V_\xi}\to K({\bf Z_d},1)\}$

The above diagram induces the following diagram.

$$
\begin{array}{ccc}
K({\bf Z},1)&\stackrel{\nu'}\rightarrow&K({\bf Z_d},1)\\ 
\uparrow_\zeta &\nearrow &\\
\hat{V_\xi}&&
\end{array}
$$

The above diagram induces the following homomorphism. 

$$
\begin{array}{cccc}
&\Omega^3(K({\bf Z},1))&\stackrel{\nu'}\rightarrow&\Omega^{3}(K({\bf Z_d},1))\\
&\rotatebox[origin=c]{90}{$\in $}& &\rotatebox[origin=c]{90}{$\in $}\\
&[(\hat{V_\xi}, \zeta)] & \mapsto& [(\hat{V_\xi}, \nu'\circ\zeta)]  
\end{array}
$$

By bordism theory 
$\Omega^3(K({\bf Z},1))\cong\Omega^3(S^1)\cong0$. 
Therefore  
$[(\hat{V_\xi},\nu'\circ\zeta)]=0$
Therefore 
$(\hat{V_\xi},\nu'\circ\zeta)=\partial(W,\tau)$ 
for a compact oriented 4-manifold $W$.   
Hence
$\widetilde\eta(\hat{V_\xi},\nu'\circ\zeta )=
\frac{1}{1}(\overline\sigma(W, \tau)-\sigma(W))$. 
Hence $\widetilde\eta(\hat{V_\xi},\nu'\circ\zeta )\in{\bf Z}$. 
Hence 
$\widetilde\eta(K,\nu)=\pi(\widetilde\eta(\hat{V_\xi},\nu'\circ\zeta ))=0$.  
Therefore Theorem 2.2 holds.

\np 
 
\footnotesize{

 }



 
  




 




\np
{\large {\bf Appendix} }

{\normalsize 
The author gives a proof of the following proposition.  

\vskip3mm
\noindent
{\bf Proposition.} 
{\it 
Let $X$ be an oriented compact $(m+2)$-dimensional manifold. 
Let $\partial X\neq\phi$. 
Let $M$ be an oriented closed $m$-dimensional manifold 
which is embedded in $X$. 
Let $M\cap\partial X\neq\phi$. 
Let $[M]=0\in H_m(X; {\bf Z})$. 
Then there is an oriented compact $(m+1)$-dimensional manifold $P$ 
such that $P$ is embedded in $X$ and that $\partial P=X$. 
}
\vskip3mm

\noindent{\bf Proof.} 
Let $\nu$ be the normal bundle of $M$ in $X$.  
By Theorem 2 in P.49 of 
\cite{Ki}   
$\nu$ is a product bundle. 
By using $\nu$ and the collar neighborhood of $\partial X$ in $X$, 
we can take a compact oriented $(m+2)$-manifold $N\subset X$ 
with the following properties. 

(1)
 $N\cong M\times D^2$. (Hence $\partial N=M\times S^1$.)

(2) 
$N\cap\partial X$
$=(\partial N)\cap(\partial X)$
$=M\cap\partial X$. 
(Hence  (Int$N)\cap\partial X=\phi$. )

Take $X-($Int$N)$. (Note $X-($Int$N)\supset \partial X$.)
There is a cell decomposition: 

$X-({\mathrm{Int}}N)$

$=(\partial N)\cup(\partial X)\cup$
(1-cells $e^1)\cup$  
(2-cells $e^2)\cup$ 
(3-cells $e^3)\cup$ 
(one 4-cell $e^4)$.

We can suppose that this decomposition has only one 0-cell $e^0$ 
which is in $(\partial N)\cap(\partial X)$.


There is a continuous map 
$s_0:(\partial N)\cup(\partial X)\to S^1$ 
with the following properties, 
where $p$ is a point in $S^1$.  

(1)  $s_0(\partial X)=p$. 
    (Hence $s_0((\partial N)\cap(\partial X))=p$ and 
    $s_0(e^0)=p$. )

(2) $s_0\vert_{\partial N}: M\times S^1\to S^1$ is a projection map 
     $(x,y)\mapsto y$.

Let $S_F^1$ be a fiber of the $S^1$-fiber bundle $\partial N=M\times S^1$. 
Since $[M]=0\in H_m(X; {\bf Z})$, 
$[S_F^1]$ generates ${\bf Z}\subset$ $H_1(X-$Int $N, \partial X;{\bf Z})$. 
(We can prove as in the proof of Theorem 3 in P.50 of 
\cite{Ki}   
) 

Let $f:H_1(X-\mathrm{Int}N, \partial X;{\bf Z})
\to H_1(X-\mathrm{Int}N, \partial X;{\bf Z})/$Tor 
be the natural projection map. 
Let $\{f([S^1_F]), u_1,...,u_k\}$ be a set of basis of 

\noindent $H_1(X-\mathrm{Int}N, \partial X;{\bf Z})/$Tor. 
Take a continuous map 

$s_1: (\partial N)\cup(\partial X)\cup$(1-cells $e^1)\to S^1$  

with the following properties. 

(1) $s_1\vert_
{(\partial N)\cup(\partial X)}=s_0$  

(2) $s_1\vert_{e^0\cup e^1}:e^0\cup e^1\to S^1$ 
satisfies the following condition: 
If $f([e^0\cup e^1])=$
\newline
$n_0\cdot f([S^1_f])+\Sigma_{j=1}^{k}n_j\cdot u_j$
$\in H_1(X-\mathrm{Int}N, \partial X;{\bf Z})/$Tor
($n_*\in{\bf Z}$), 
\newline  
then deg$(s_1\vert_{e^0\cup e^1})=n_0$.

Note that, if a circle $C$ is nul-homologous in 
$(\partial N)\cup(\partial X)\cup$(1-cells $e^1)$, 
then deg$(s_1\vert_C)=0$. 

\noindent
{\bf Claim.}  
{\it 
There is a continuous map 

$s_2:(\partial N)\cup(\partial X)\cup$(1-cells $e^1)\cup$(2-cells $e^2)
\to S^1$  

such that 
$s_2\vert_{(\partial N)\cup(\partial X)\cup({\mathrm{1-cells}}\quad e^1)}$
=$s_1$.
}

\noindent
{\bf Proof.} 
 It is trivial that $[\partial e^2]=0$
$\in 
H_1((\partial N)\cup(\partial X)\cup$(1-cells $e^1);{\bf Z})$. 
Hence deg($s_1\vert_{\partial e^2}$) is zero.  
Hence $s_1\vert_{\partial e^2}$ extends to $e^2$. 
Hence the above Claim holds. 

The continuous map $s_2$ extends to a continuous map 
$s:X-({\mathrm{Int}}N)\to S^1$ 
since $\pi_l(S^1)=0 (l\geq2)$. 
We can suppose $s$ is a smooth map.

Let $q\neq p$. Let $q$ be a regular value. 
Hence $s^{-1}(q)$ be an oriented compact manifold. 
$\partial\{s^{-1}(q)\}\subset\{(\partial N)\cup\partial X\}$.   
Since $q\neq p$, $s^{-1}(q)\cap\partial X=\phi.$  
Hence $\partial\{s^{-1}(q)\}\subset\partial N$.
Furthermore we have 
$s^{-1}(q)\cap \partial N=\partial \{s^{-1}(q)\}  
=M\times \{r\}$, where $r$ is a point in $S^1$. 
By using $N$ and $s^{-1}(q)$, Proposition holds. 

}

\np

\pagestyle{empty}

\unitlength 0.1in
\begin{picture}(56.10,47.00)(8.50,-47.60)
%
\special{pn 8}%
\special{ar 3510 320 560 250  0.0000000 6.2831853}%
%
\special{pn 20}%
\special{pa 3660 3950}%
\special{pa 3660 2730}%
\special{fp}%
%
\special{pn 20}%
\special{pa 3340 3970}%
\special{pa 3340 2740}%
\special{fp}%
%
\special{pn 20}%
\special{pa 3350 4000}%
\special{pa 3360 3970}%
\special{pa 3384 3949}%
\special{pa 3413 3935}%
\special{pa 3444 3926}%
\special{pa 3475 3921}%
\special{pa 3507 3920}%
\special{pa 3539 3923}%
\special{pa 3570 3930}%
\special{pa 3601 3940}%
\special{pa 3627 3958}%
\special{pa 3647 3983}%
\special{pa 3648 4014}%
\special{pa 3630 4040}%
\special{pa 3603 4058}%
\special{pa 3573 4070}%
\special{pa 3542 4077}%
\special{pa 3510 4080}%
\special{pa 3478 4079}%
\special{pa 3447 4075}%
\special{pa 3416 4066}%
\special{pa 3387 4053}%
\special{pa 3362 4032}%
\special{pa 3350 4003}%
\special{pa 3350 4000}%
\special{sp}%
%
\special{pn 20}%
\special{pa 3340 360}%
\special{pa 3340 1580}%
\special{fp}%
%
\special{pn 20}%
\special{pa 3660 340}%
\special{pa 3660 1570}%
\special{fp}%
%
\special{pn 20}%
\special{ar 3500 310 150 80  0.0000000 6.2831853}%
%
\special{pn 8}%
\special{ar 3510 4080 560 250  0.0000000 6.2831853}%
%
\special{pn 8}%
\special{pa 4080 330}%
\special{pa 4080 4080}%
\special{fp}%
%
\special{pn 8}%
\special{pa 2950 340}%
\special{pa 2950 4080}%
\special{fp}%
%
\special{pn 8}%
\special{ar 1410 330 560 250  0.0000000 6.2831853}%
%
\special{pn 8}%
\special{ar 1410 4090 560 250  0.0000000 6.2831853}%
%
\special{pn 8}%
\special{pa 1980 340}%
\special{pa 1980 4090}%
\special{fp}%
%
\special{pn 8}%
\special{pa 850 350}%
\special{pa 850 4090}%
\special{fp}%
%
\special{pn 8}%
\special{ar 5890 310 560 250  0.0000000 6.2831853}%
%
\special{pn 8}%
\special{ar 5890 4070 560 250  0.0000000 6.2831853}%
%
\special{pn 8}%
\special{pa 6460 320}%
\special{pa 6460 4070}%
\special{fp}%
%
\special{pn 8}%
\special{pa 5330 330}%
\special{pa 5330 4070}%
\special{fp}%
%
\special{pn 20}%
\special{ar 3500 1590 150 80  0.0000000 6.2831853}%
%
\special{pn 20}%
\special{ar 3500 2710 150 80  0.0000000 6.2831853}%
%
\special{pn 20}%
\special{ar 3520 2070 560 250  0.0000000 6.2831853}%
%
\special{pn 8}%
\special{pa 3660 1600}%
\special{pa 5730 1600}%
\special{dt 0.045}%
\special{pa 5730 1600}%
\special{pa 5729 1600}%
\special{dt 0.045}%
%
\special{pn 8}%
\special{pa 3680 2720}%
\special{pa 5750 2720}%
\special{dt 0.045}%
\special{pa 5750 2720}%
\special{pa 5749 2720}%
\special{dt 0.045}%
%
\special{pn 20}%
\special{ar 5850 1600 150 80  0.0000000 6.2831853}%
%
\special{pn 20}%
\special{ar 5850 2730 150 80  0.0000000 6.2831853}%
%
\special{pn 20}%
\special{pa 5690 1610}%
\special{pa 5690 2700}%
\special{fp}%
%
\special{pn 20}%
\special{pa 6020 1640}%
\special{pa 6020 2730}%
\special{fp}%
\put(12.0000,-46.0000){\makebox(0,0)[lb]{t=-0.5}}%
\put(32.0000,-46.0000){\makebox(0,0)[lb]{t=0}}%
\put(56.0000,-46.0000){\makebox(0,0)[lb]{t=0.5}}%
\put(30.3000,-49.3000){\makebox(0,0)[lb]{Figure 1.1}}%
%
\special{pn 8}%
\special{ar 3510 320 560 250  0.0000000 6.2831853}%
%
\special{pn 8}%
\special{pa 2950 340}%
\special{pa 2950 4080}%
\special{fp}%
%
\special{pn 8}%
\special{pa 4080 330}%
\special{pa 4080 4080}%
\special{fp}%
%
\special{pn 8}%
\special{ar 3510 4080 560 250  0.0000000 6.2831853}%
%
\special{pn 20}%
\special{ar 3500 310 150 80  0.0000000 6.2831853}%
%
\special{pn 20}%
\special{pa 3660 340}%
\special{pa 3660 1570}%
\special{fp}%
%
\special{pn 20}%
\special{pa 3340 360}%
\special{pa 3340 1580}%
\special{fp}%
%
\special{pn 20}%
\special{pa 3350 4000}%
\special{pa 3360 3970}%
\special{pa 3384 3949}%
\special{pa 3413 3935}%
\special{pa 3444 3926}%
\special{pa 3475 3921}%
\special{pa 3507 3920}%
\special{pa 3539 3923}%
\special{pa 3570 3930}%
\special{pa 3601 3940}%
\special{pa 3627 3958}%
\special{pa 3647 3983}%
\special{pa 3648 4014}%
\special{pa 3630 4040}%
\special{pa 3603 4058}%
\special{pa 3573 4070}%
\special{pa 3542 4077}%
\special{pa 3510 4080}%
\special{pa 3478 4079}%
\special{pa 3447 4075}%
\special{pa 3416 4066}%
\special{pa 3387 4053}%
\special{pa 3362 4032}%
\special{pa 3350 4003}%
\special{pa 3350 4000}%
\special{sp}%
%
\special{pn 20}%
\special{pa 3340 3970}%
\special{pa 3340 2740}%
\special{fp}%
%
\special{pn 20}%
\special{pa 3660 3950}%
\special{pa 3660 2730}%
\special{fp}%
%
\special{pn 8}%
\special{ar 1410 330 560 250  0.0000000 6.2831853}%
%
\special{pn 8}%
\special{ar 1410 4090 560 250  0.0000000 6.2831853}%
%
\special{pn 8}%
\special{pa 1980 340}%
\special{pa 1980 4090}%
\special{fp}%
%
\special{pn 8}%
\special{pa 850 350}%
\special{pa 850 4090}%
\special{fp}%
%
\special{pn 8}%
\special{ar 5890 310 560 250  0.0000000 6.2831853}%
%
\special{pn 8}%
\special{ar 5890 4070 560 250  0.0000000 6.2831853}%
%
\special{pn 8}%
\special{pa 6460 320}%
\special{pa 6460 4070}%
\special{fp}%
%
\special{pn 8}%
\special{pa 5330 330}%
\special{pa 5330 4070}%
\special{fp}%
\end{picture}%

\np
\unitlength 0.1in
\begin{picture}(56.10,47.10)(8.50,-47.70)
%
\special{pn 8}%
\special{ar 3510 320 560 250  0.0000000 6.2831853}%
%
\special{pn 20}%
\special{pa 3660 3950}%
\special{pa 3660 2730}%
\special{fp}%
%
\special{pn 20}%
\special{pa 3340 3970}%
\special{pa 3340 2740}%
\special{fp}%
%
\special{pn 20}%
\special{pa 3350 4000}%
\special{pa 3360 3970}%
\special{pa 3384 3949}%
\special{pa 3413 3935}%
\special{pa 3444 3926}%
\special{pa 3475 3921}%
\special{pa 3507 3920}%
\special{pa 3539 3923}%
\special{pa 3570 3930}%
\special{pa 3601 3940}%
\special{pa 3627 3958}%
\special{pa 3647 3983}%
\special{pa 3648 4014}%
\special{pa 3630 4040}%
\special{pa 3603 4058}%
\special{pa 3573 4070}%
\special{pa 3542 4077}%
\special{pa 3510 4080}%
\special{pa 3478 4079}%
\special{pa 3447 4075}%
\special{pa 3416 4066}%
\special{pa 3387 4053}%
\special{pa 3362 4032}%
\special{pa 3350 4003}%
\special{pa 3350 4000}%
\special{sp}%
%
\special{pn 20}%
\special{pa 3340 360}%
\special{pa 3340 1580}%
\special{fp}%
%
\special{pn 20}%
\special{pa 3660 340}%
\special{pa 3660 1570}%
\special{fp}%
%
\special{pn 20}%
\special{ar 3500 310 150 80  0.0000000 6.2831853}%
%
\special{pn 8}%
\special{ar 3510 4080 560 250  0.0000000 6.2831853}%
%
\special{pn 8}%
\special{pa 4080 330}%
\special{pa 4080 4080}%
\special{fp}%
%
\special{pn 8}%
\special{pa 2950 340}%
\special{pa 2950 4080}%
\special{fp}%
%
\special{pn 8}%
\special{ar 1410 330 560 250  0.0000000 6.2831853}%
%
\special{pn 8}%
\special{ar 1410 4090 560 250  0.0000000 6.2831853}%
%
\special{pn 8}%
\special{pa 1980 340}%
\special{pa 1980 4090}%
\special{fp}%
%
\special{pn 8}%
\special{pa 850 350}%
\special{pa 850 4090}%
\special{fp}%
%
\special{pn 8}%
\special{ar 5890 310 560 250  0.0000000 6.2831853}%
%
\special{pn 8}%
\special{ar 5890 4070 560 250  0.0000000 6.2831853}%
%
\special{pn 8}%
\special{pa 6460 320}%
\special{pa 6460 4070}%
\special{fp}%
%
\special{pn 8}%
\special{pa 5330 330}%
\special{pa 5330 4070}%
\special{fp}%
%
\special{pn 20}%
\special{ar 3500 1590 150 80  0.0000000 6.2831853}%
%
\special{pn 20}%
\special{ar 3500 2710 150 80  0.0000000 6.2831853}%
%
\special{pn 20}%
\special{ar 3520 2070 560 250  0.0000000 6.2831853}%
\put(12.0000,-46.0000){\makebox(0,0)[lb]{t=-0.5}}%
\put(32.0000,-46.0000){\makebox(0,0)[lb]{t=0}}%
\put(56.0000,-46.0000){\makebox(0,0)[lb]{t=0.5}}%
%
\special{pn 8}%
\special{ar 3510 320 560 250  0.0000000 6.2831853}%
%
\special{pn 8}%
\special{pa 2950 340}%
\special{pa 2950 4080}%
\special{fp}%
%
\special{pn 8}%
\special{pa 4080 330}%
\special{pa 4080 4080}%
\special{fp}%
%
\special{pn 8}%
\special{ar 3510 4080 560 250  0.0000000 6.2831853}%
%
\special{pn 20}%
\special{ar 3500 310 150 80  0.0000000 6.2831853}%
%
\special{pn 20}%
\special{pa 3660 340}%
\special{pa 3660 1570}%
\special{fp}%
%
\special{pn 20}%
\special{pa 3340 360}%
\special{pa 3340 1580}%
\special{fp}%
%
\special{pn 20}%
\special{pa 3350 4000}%
\special{pa 3360 3970}%
\special{pa 3384 3949}%
\special{pa 3413 3935}%
\special{pa 3444 3926}%
\special{pa 3475 3921}%
\special{pa 3507 3920}%
\special{pa 3539 3923}%
\special{pa 3570 3930}%
\special{pa 3601 3940}%
\special{pa 3627 3958}%
\special{pa 3647 3983}%
\special{pa 3648 4014}%
\special{pa 3630 4040}%
\special{pa 3603 4058}%
\special{pa 3573 4070}%
\special{pa 3542 4077}%
\special{pa 3510 4080}%
\special{pa 3478 4079}%
\special{pa 3447 4075}%
\special{pa 3416 4066}%
\special{pa 3387 4053}%
\special{pa 3362 4032}%
\special{pa 3350 4003}%
\special{pa 3350 4000}%
\special{sp}%
%
\special{pn 20}%
\special{pa 3340 3970}%
\special{pa 3340 2740}%
\special{fp}%
%
\special{pn 20}%
\special{pa 3660 3950}%
\special{pa 3660 2730}%
\special{fp}%
%
\special{pn 8}%
\special{ar 1410 330 560 250  0.0000000 6.2831853}%
%
\special{pn 8}%
\special{ar 1410 4090 560 250  0.0000000 6.2831853}%
%
\special{pn 8}%
\special{pa 1980 340}%
\special{pa 1980 4090}%
\special{fp}%
%
\special{pn 8}%
\special{pa 850 350}%
\special{pa 850 4090}%
\special{fp}%
%
\special{pn 8}%
\special{ar 5890 310 560 250  0.0000000 6.2831853}%
%
\special{pn 8}%
\special{ar 5890 4070 560 250  0.0000000 6.2831853}%
%
\special{pn 8}%
\special{pa 6460 320}%
\special{pa 6460 4070}%
\special{fp}%
%
\special{pn 8}%
\special{pa 5330 330}%
\special{pa 5330 4070}%
\special{fp}%
%
\special{pn 8}%
\special{pa 3320 1590}%
\special{pa 1590 1590}%
\special{dt 0.045}%
\special{pa 1590 1590}%
\special{pa 1591 1590}%
\special{dt 0.045}%
%
\special{pn 8}%
\special{pa 3330 2740}%
\special{pa 1620 2740}%
\special{dt 0.045}%
\special{pa 1620 2740}%
\special{pa 1621 2740}%
\special{dt 0.045}%
%
\special{pn 20}%
\special{ar 1420 1590 150 80  0.0000000 6.2831853}%
%
\special{pn 20}%
\special{ar 1420 2720 150 80  0.0000000 6.2831853}%
%
\special{pn 20}%
\special{pa 1260 1600}%
\special{pa 1260 2690}%
\special{fp}%
%
\special{pn 20}%
\special{pa 1590 1630}%
\special{pa 1590 2720}%
\special{fp}%
\put(30.0000,-49.4000){\makebox(0,0)[lb]{Figure 1.2}}%
\end{picture}%

\np
\input 6.1.tex

\np
\input 6.2.tex

\np
\input 6.3.tex

\np
\input 6.4.tex

\np
\input 6.5.tex

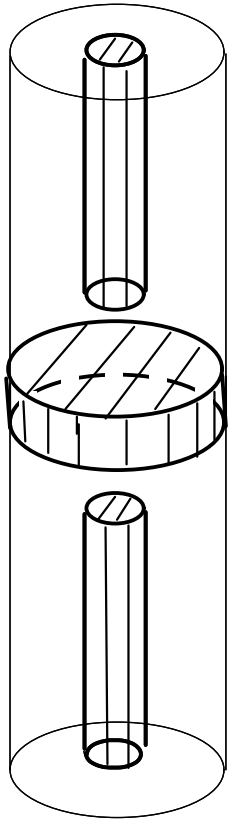
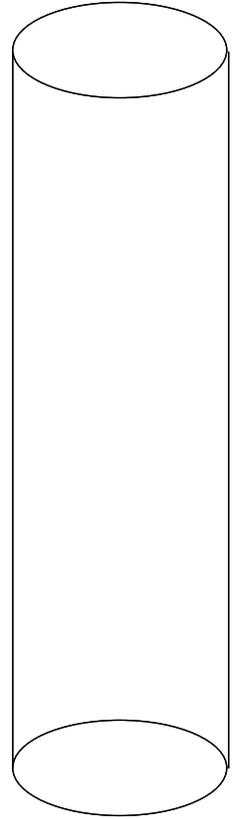
\end{document}